\documentclass[11pt]{article}
\usepackage{amsmath,amsfonts,latexsym,amssymb,amscd, graphicx,pictexwd, mathrsfs}
\newcommand{\w}{{\mathrm w}}
\newcommand{\s}{{\mathrm s}}
\newcommand{\f}{{\mathrm f}}
\newcommand{\rb}{{\mathrm b}}
\newcommand{\m}{{\mathrm m}}
\newcommand{\ra}{{\mathrm a}}
\newcommand{\rt}{{\mathrm t}}

\renewcommand{\P}{{\mathbb P}}
\newcommand{\E}{{\mathbb E}\,}

\newcommand{\R}{{\mathbb R}}

\newcommand{\ZZ}{{\mathbb Z}}

\newcommand{\cc}{{\mathbf c}}
\newcommand{\zz}{{\mathbf z}}

\newcommand{\xx}{{\mathbf x}}
\newcommand{\yy}{{\mathbf y}}
\newcommand{\0}{{\mathbf 0}}

\newcommand{\argmax}{\mathop{\rm arg\,max}}

\newcommand{\Exp}{{\rm Exp}}

\newcommand{\Var}{{\rm Var\,}}

\newtheorem{thm}{Theorem}

\newtheorem{coro}{Corollary}[section]
\newtheorem{lem}{Lemma}[section]

\newtheorem{rem}{Remark}

\newtheorem{ass}{Assumption}
\numberwithin{equation}{section}

\begin{document}
\title{Local Behavior of Airy Processes}
\author{Leandro P. R. Pimentel}
\date{\today}
\maketitle

\begin{abstract}  
The Airy processes describe spatial fluctuations in wide range of growth models, where each particular Airy process arising in each case depends on the geometry of the initial profile. We show how the coupling method, developed in the last-passage percolation context, can be used to prove that several types of Airy processes have a continuous version, and behave locally like a Brownian motion. We further extend these results to an Airy sheet, by proving existence of a continuous version and local convergence to additive Brownian motion.  
\end{abstract}

\section{Introduction and Main Results}
The Airy processes describe spatial fluctuations in wide range of growth models, where each particular Airy process arising in each case depends on the geometry of the initial profile. These processes are at the centre of the KPZ universality class, and it is conjectured to describe the fluctuations of wide class of interface growth models characterized by satisfying a local slope dependent growth rate and a smoothing mechanism, combined with space-time random forcing with rapid decay of correlations. The KPZ equation \cite{AmCoQu,KaPaZh} is a canonical example of such model, giving its name to the universality class. The TASEP height process and the last-passage percolation model (exponential or geometric) \cite{Jo} are also examples of interface growth processes that present KPZ fluctuations.          
\newline

In \cite{CaPi}, the coupling method was applied to study local fluctuations of point-to-point last-passage percolation (LPP) times and its scaling limit, the Airy$_2$ process. The main technique relies on a local comparison lemma that allows us to sandwich local fluctuations of point-to-point LPP times in between Brownian local fluctuations of the equilibrium regime. The local comparison lemma can be seen as a property involving a basic coupling of LPP models started from narrow wedge (point-to-point) and equilibrium initial profiles. The main contribution of this paper is to use an extension of the local comparison property for two arbitrary initial profiles, to show existence of continuous versions and Brownian local fluctuations for several types of Airy processes, where each particular Airy process arising in each case depends on the geometry of initial profile (Theorem \ref{main}). In addition, we state that Brownian fluctuations of local increments occur under sub-KPZ scaling (Remark \ref{subKPZ}), and give a brief script for the proof. The coupling method will be further applied to prove existence of continuous versions of a two two parameter fluctuation process, called Airy sheet \cite{MaQuRe}, and local convergence to a sum of two independent Brownian motions  (Theorem \ref{mainSheet}) \footnote{The coupling method is also an useful  tool to study the long time behavior of the KPZ fixed point and prove ergodicity  \cite{MaQuRe,Pi} of its local increments.}.
\newline

In the course of this article will need to verify a crucial property  (Assumption \ref{exitcontrol}), which concerns localization of a LPP maximizer under the KPZ scaling. This property was already studied in \cite{ChFeSp,CoLiWa,Jo}, to prove variational formulas for Airy processes, by using asymptotic analysis of a kernel that arises from exact formulas for LPP models with geometric or exponential weights. We will follow the coupling method approach to prove localization for some examples of initial profiles (Theorem \ref{local}), and it will parallel the ideas in \cite{BaCaSe,CaGr} developed for the equilibrium situation.   
\newline

\subsection{Local fluctuations of the LPP model}
We start with some definitions. For $\xx,\yy\in\ZZ^2$ we define that $\xx \leq \yy$, for $\xx=(i,j)\in\ZZ^2$ and $\yy=(k,l)$, if $i\leq k$ and $j\leq l$. We say that a sequence $\pi=(\xx_0,\xx_1,\cdots,x_n)$ is an up-right path from $\xx$ to $\yy$, with $\xx\leq \yy$, if $\xx_{m+1}-\xx_{m}\in\{(1,0),(0,1)\}$, $\xx_0=\xx$ and $\xx_n=\yy$. We denote $\Pi(\xx,\yy)$ the set of all up-right paths from $\xx$ to $\yy$.  The random environment in our setting is given by a collection $\omega\equiv\left\{\omega_{i,j}\,:\,i+j>0\right\}$ of i.i.d. random variables (passage times) with exponential  distribution of parameter $1$. Given $\pi=(\xx_0,\xx_1,\cdots,\xx_n)\in\Pi(\xx,\yy)$ the passage time is set as 
$$\omega(\pi):= \sum_{i=1}^n \omega_i\,.$$
Notice that we do not include the passage time at $\xx_0=\xx$. For $\xx\leq \yy$, with $\xx=(i,j)$ and  $i+j>0$, the point-to-point last-passage percolation time is defined as
$$L(\xx,\yy):=\max_{\pi\in\Pi(\xx,\yy)}\omega(\pi)\,.$$
  
Denote $L(\xx)\equiv L(\0,\xx)$, $L_k(\xx)\equiv L((k,-k),\xx)$ and set 
$$C^\xx:=\left\{k\in\ZZ\,:\, (k,-k)\leq \xx\right\}\,.$$
Given a profile 
$$\rb:\ZZ\to\ZZ\cup\{-\infty\}\mbox{ with }\rb(0)=0\,,$$ 
the last-passage percolation time with initial boundary profile $\rb$ is defined as 
$$L^{\rb}(\xx):=\max_{k\in C^\xx} \left\{ \rb(k) + L_k(\xx)\right\}\,.$$
The last-passage percolation model generates a discrete time Markov process $(M^\rb_n)_{n\geq 0}$ defined as 
$$M^\rb_n(k):= L^\rb[k]_n-L^\rb[0]_n\,,\,\mbox { for }\,k\in\ZZ\,.$$
(Note that $M^\rb_0\equiv\rb$.) The increments $L^\rb[k]_n-L^\rb[0]_n$ are define along $[k]_n=(n+k,n-k)$, so that $n$ represents the time parameter. For a real number $x\in [-n,n]$ we denote 
$$[x]_n\equiv [\lfloor x\rfloor]_n\equiv (n+\lfloor x\rfloor , n-\lfloor x\rfloor)\,.$$  
For a fixed constant $C>0$, define the process 
$$\Delta_n^\rb(u):=\frac{L^\rb[un^{2/3}]_n-L^\rb[0]_n}{2^{3/2}n^{1/3}}\,,u\in[0,C]\,$$ 
(where $n\geq C^3$). We could have defined a continuous one by linearly interpolating the values, but this would not be relevant for our purposes. Our goal is to prove tightness of $\Delta_n^{\rb}$ and a local functional convergence of any weak limit to Brownian motion. 
\newline

Let $\Exp_1(1-\rho)$ and $\Exp_2(\rho)$ be independent random variables with exponential distributions of parameter $1-\rho$ and $\rho$, respectively, and define 
$$X(\rho)=\Exp_1(1-\rho)-\Exp_2(\rho)\,.$$
The unique family of time stationary measures (with ergodic space increments) for the LPP Markov process is the one induced by a sum of i.i.d. copies $\zeta_k$ of $X(\rho)$. If one sets $\s_\rho(0)=0$ and 
\begin{equation}\label{defstat}
\s_\rho(k)=\left\{\begin{array}{lll}
\sum_{i=k+1}^{0}-\zeta_i & \mbox{ for } k<0\\
0 & \mbox{ for } k=0\\
\sum_{i=1}^k\zeta_i &\mbox{ for } k>0\,,\end{array}\right.
\end{equation} 
then \cite{BaCaSe}
\begin{equation}\label{stat}
M^{\s_\rho}_n(\cdot)\stackrel{dist.}{=}\s_\rho(\cdot)\,,\,\forall\,n\geq 0\,.
\end{equation}
As a corollary of the functional central limit theorem for sums of i.i.d. random variables, in the stationary regime with $\rho=1/2$, we have that 
\begin{equation}\label{ConvEqui}
\lim_{n\to\infty}\Delta_n^{\s_{1/2}}(\cdot)\stackrel{dist.}{=} B(\cdot)\,,
\end{equation}
in the Skorohod topology of cadlag functions on compact sets, where $B$ is a standard Brownian motion. As we shall see next, we will be able to study $\Delta_n^\rb$ by comparing its local behavior with the stationary regime. 
\newline

To state the main result of this paper we need to introduce the exit-point location,  
$$Z^{\rb}(\xx) :=\max\arg\max_{k\in C_{\xx}} \left\{ \rb(k) + L_k(\xx)\right\}\,,$$
in such a way that, 
$$L^{\rb}(\xx)=\rb\left(Z^{\rb}(\xx)\right)+L_{Z^{\rb}(\xx)}(\xx)\,.$$
The location of the maximizer might be not unique. For instance, if one sets $\rb\equiv 0$ (flat), then the geodesic to $(1,0)$ can either do $(0,0)\to(1,0)$ or $(0,0)\to(1,-1)\to(1,0)$, since the weights along the boundary are the same. In this example, the exit point is $1$ (the right most). The exit point location is the key to compare the evolution of $M_n^\rb$ with the equilibrium regime $M_n^{\s_{1/2}}$ (see Lemma \ref{comparison}). To keep the evolution close enough to equilibrium, we will need the following assumption.   

\begin{ass}\label{exitcontrol}   
Let $C\geq 0$ and for $r\geq 0$ define
$$R_C(r):= \limsup_{n\to\infty}\P\left(|Z^\rb[C n^{2/3}]_n|\geq rn^{2/3} \right)\,.$$
Then
$$\limsup_{r\to\infty}R_C(r)=0\,,\,\mbox{for all $C\geq 0$}\,.$$
\end{ass} 

\begin{thm}\label{main}
Under Assumption \ref{exitcontrol}, the collection $\left\{\Delta^\rb_{n}(\cdot)\,:\,n\geq 1\right\}$ is relatively compact in the Skorohod topology of cadlag functions on compact sets,  and any weak limit point is continuous almost surely. For a weak limit point $\Delta^\rb(\cdot)$ of $\left\{\Delta^\rb_{n}(\cdot)\,:\,n\geq 1\right\}$ define $\Delta^{\rb,\epsilon}(x):=\epsilon^{-1/2}\Delta^\rb(\epsilon x)$. Then       
\begin{equation}\label{loconv}
\lim_{\epsilon\downarrow 0}\Delta^{\rb,\epsilon}(\cdot)\stackrel{dist.}{=} B(\cdot)\,,
\end{equation}
in the topology of continuous functions on compact sets, where $(B(x)\,,x\geq 0)$ is a standard Brownian motion.
\end{thm}

\begin{rem}\label{subKPZ}
Local fluctuations under sub-KPZ scaling is defined by fixing $\gamma\in(0,2/3)$ and setting
$$\Delta_{\gamma,n}^{\rb}(u):=\frac{L^\rb[un^{\gamma}]_n-L^\rb[0]_n}{2^{3/2}n^{\gamma/2}}\,,u\in[0,C]\,.$$ 
As in \cite{CaPi} (see Theorems 2 and 5 there), one can see that the script to prove Theorem \ref{main} can be adapted to show that, under Assumption \ref{exitcontrol},
\begin{equation}\label{subloconv}
\lim_{n\to\infty}\Delta_{\gamma,n}^{\rb}(\cdot)\stackrel{dist.}{=} B(\cdot)\,,
\end{equation}
in the topology of in the Skorohod topology of cadlag functions on compact sets.
\end{rem}

\subsection{Examples of initial profiles}

Assumption \ref{exitcontrol} concerns the localization of the exit point in the $n^{2/3}$ scale  around the origin. This behavior is crucial when one is proving variational formulas for Airy processes \cite{Jo}. In the last-passage percolation model with geometric weights, and random walk type of initial profile, the analog of Assumption \ref{exitcontrol} can be obtained from (98) and (99) in \cite{CoLiWa}. For exponential weights, one can use the estimates in the proof of Theorem 2.2 in \cite{ChFeSp} to get localization for a large range of random walk type of initial profiles (see Remark 2.7 there). A common ingredient in \cite{ChFeSp,CoLiWa,Jo} is the use of asymptotic analysis of a kernel that arises from exact formulas for LPP models with geometric or exponential weights. In this article we will use a different approach to prove localization that parallels the coupling method developed in \cite{BaCaSe,CaGr}. See Section \ref{check} for the proofs.  

\begin{thm}\label{local}
Assumption \ref{exitcontrol} holds for all profiles below.
\end{thm}             
 
\subsubsection{Narrow Wedge Profile}
If
$$\w(k)=\left\{\begin{array}{ll}0 &\mbox{ for } k=0\\
-\infty& \mbox{ for } k\neq 0\,,\end{array}\right.$$ 
then 
$$L^{\w}(\xx)=L(\xx)\mbox{ and }Z^\w(\xx)\equiv 0\,.$$
Define 
$$H^\w_{n}(u) = \frac{L^\w[2^{2/3} un^{2/3}]_n-4 n}{2^{4/3} n^{1/3}}\,,u\in[0,C]\,.$$
Then 
$$H^{\w}_{n}(u)=H^{\w}_{n}(0)+2^{1/6}\Delta^\w_n(2^{2/3} u)\,.$$ The Airy$_2$ process arise as the limit of $A^\w_{n}$ \cite{BoFe,PrSp}:
$$\lim_{n\to\infty}H^\w_{n}(u) \stackrel{dist.}{=} A_2(u)-u^2\,,$$
Clearly Assumption \ref{exitcontrol} is true in this case, and we can use Theorem \ref{main} to show tightness of $H^\w_n$ and local (functional) convergence of the Airy$_2$ process to Brownian Motion:
\begin{equation}\label{locA2}
\lim_{\epsilon\downarrow 0}\epsilon^{-1/2}\left(A_2(\epsilon x)-A_2(0)\right)\stackrel{dist.}{=} \sqrt{2}B(x)\,.
\end{equation}
This was actually the main result in \cite{CaPi}, in a slightly different context, where the local comparison was use for the first time. Different approaches to prove local convergence were developed in \cite{CoHa,Ha}. For any initial profiles that produces a rarefaction fan (wedge type of profile), the exit point will stabilize close to origin \cite{CaPi1}. This can be used to prove Assumption \ref{exitcontrol} and, as a consequence, to apply Theorem \ref{main} in this context as well.

\subsubsection{Flat Profile}
The flat profile is defined as $\f(k)=0$ for all $k\in\ZZ$ (line to point last-passage percolation). Define 
$$H^\f_{n}(u) = \frac{L^\f[2^{2/3} un^{2/3}]_n-\mu n}{2^{4/3} n^{1/3}}\,,u\in[0,C]\,.$$
Then
$$H^{\f}_{n}(u)=H^{\f}_{n}(0)+2^{1/6}\Delta^\f_n(2^{2/3} u)\,.$$ The Airy$_1$ process arise as the limit of $H^\f_{n}$:
$$\lim_{n\to\infty}H^\f_{n}(u) \stackrel{dist.}{=} A_1(u)\,,$$
in the sense of convergence of finite-dimensional distributions \cite{BoFePrSa}. 
\newline

Theorem \ref{main}, together with Theorem \ref{local}, implies tightness of $H^\f_{n}$, and local (functional) convergence of the Airy$_1$ process to Brownian Motion:
\begin{equation}\label{locA1}
\lim_{\epsilon\downarrow 0}\epsilon^{-1/2}\left(A_1(\epsilon x)-A_1(0)\right)\stackrel{dist.}{=} \sqrt{2}B(x)\,.
\end{equation}
We believe that this the first work that brings tightness of $H^\f_{n}$ and functional convergence of local fluctuations of $A_1$. Finite dimensional convergence of the Airy$_1$ was first prove in \cite{QuRe}. We note that 
$$\P\left(\max_{k\in [-rn^{2/3},rn^{2/3}]} \left\{L_k[0]_n\right\}\neq L^\f[0]_n\right)=\P\left(|Z^\f[0]_n|>rn^{2/3}\right)\,.$$
Therefore, Theorem \ref{local} implies that 
\begin{equation}\label{exitdist}
\lim_{n\to\infty}\frac{Z^\f[0]_n}{2^{2/3}n^{2/3}}\stackrel{dist.}{=} \argmax_{u}\left\{A_2(u)-u^2\right\}\,.
\end{equation}
See Proposition 1.4 and Theorem 1.6 in \cite{Jo} for analogous results for LPP with geometric weights.

\subsubsection{Mixed Profiles}
Mixed initial profiles can be obtained by placing one condition on each half of $\ZZ$: 
$$\w\f(k)=\left\{\begin{array}{ll}-\infty & \mbox{ for } k<0\\
0&\mbox{ for } k\geq 0\,,\end{array}\right.\mbox{\,(wedge to flat)}$$
$$\w\s(k)=\left\{\begin{array}{lll} -\infty & \mbox{ for } k<0\\
0 &\mbox{ for } k=0\\
\sum_{i=1}^k\zeta_i &\mbox{ for } k>0\,,\end{array}\right.\mbox{\,(wedge to stationary)}$$ 
$$\f\s(k)=\left\{\begin{array}{ll} 0 & \mbox{ for } k\leq 0\\
\sum_{i=1}^k\zeta_i&\mbox{ for } k>0\,.\end{array}\right.\mbox{\,(flat to stationary)}$$ 
For each case there is a specific Airy process, denoted Airy$_{2\to 1}$ (wedge to flat), Airy$_{2\to 0}$ (wedge to stationary) and Airy$_{1\to 0}$ (flat to stationary).

Theorem \ref{main}, together with Theorem \ref{local}, implies relative compactness of each rescaled processes, and local (functional) convergence to Brownian Motion of the respective Airy process. We believe that this the first work that brings functional convergence of local fluctuations for mixed profiles. We note again that Theorem \ref{local} implies that the maximum is attained in a compact set with high probability, which allow us to related the distribution of the exit point with a variational problem involving the Airy$_2$ process (as in \eqref{exitdist}).
 
\subsubsection{A remark on the Short Scale Fluctuations of Maximizers}
For fixed $t\in(0,1]$ consider the scaling operator on functions $f(x)$ given by 
$$S_{t}(f)(x):=t^{1/3} f(t^{-2/3}x)\,.$$
Let $A_2^{(1)},A^{(2)}_2$ be two independent Airy$_2$ processes and define
$$X_t:=\arg\max_{x\in\R}\left\{S_{t}\left(A_2^{(1)}\right)(x)+S_{1-t}\left(A_2^{(2)}\right)(x)-\frac{x^2}{t(1-t)} \right\}\,.$$
Using Theorem 4 \cite{Pi1}, it is not hard to see that this maximizer is a.s. unique, since the sum of two independent stationary process is a stationary process. The location $X_t$, for $t\in(0,1)$, describes the limit fluctuations, under the $n^{2/3}$ scaling, of the intersection between the geodesic from $(0,0)$ to $(n,n)$ (point-to-point) and the anti-diagonal $i+j=2\lfloor nt\rfloor$.  
\newline

If instead one is interested in the fluctuations of the intersection between the point-to-line geodesic\footnote{From $(0,0)$ to the anti-diagonal $i+j= n$.} and the anti-diagonal $i+j=2\lfloor nt\rfloor$, we get 
$$Y_t:=\arg\max_{y\in\R}\left\{S_{t}\left(A_2\right)(y)+S_{1-t}\left(A_1\right)(y)-\frac{y^2}{t} \right\}\,,$$
where $A_2$ and $A_1$ are two independent Airy$_2$ and Airy$_1$ processes, respectively. Theorem 4 \cite{Pi1} gives uniqueness of $Y_t$ by the same reason as before.
\newline

Functional local convergence of the Airy$_2$ and Airy$_1$ processes implies that $t^{-2/3}X_t$ and $t^{-2/3}Y_t$ also have a limit behavior. In view of the definition of $X_t$ and $Y_t$:
\begin{eqnarray}
\nonumber t^{-2/3}X_t&=&\arg\max_{z\in\R}\left\{t^{1/3}A_2^{(1)}(z)+(1-t)^{1/3}A_2^{(2)}\left((1-t)^{-2/3}t^{2/3}z\right)-t^{1/3}\frac{z^2}{(1-t)} \right\}\\
\label{locA21}&=&\arg\max_{z\in\R}\left\{A_2^{(1)}(z)+\epsilon^{-1/2}_t\left[A_2^{(2)}(\epsilon_t z)-A_2(0)\right]- (1+\epsilon^{3/2}_t)z^2) \right\}\,,
\end{eqnarray}
where $\epsilon_t:=(1-t)^{-2/3}t^{2/3}$, and similarly,
\begin{equation}\label{locA11}
t^{-2/3}Y_t=\arg\max_{z\in\R}\left\{A_2(z)+\epsilon_t^{-1/2}\left[A_1(\epsilon_t z)-A_1(0)\right]-z^2 \right\}\,.
\end{equation}
By \eqref{locA2} and \eqref{locA1}, one expects that   
$$\lim_{t\downarrow 0} t^{-2/3}X_t\stackrel{dist.}{=}Z\,\mbox{ and }\,\lim_{t\downarrow 0} t^{-2/3}Y_t\stackrel{dist.}{=}Z\,,$$
where 
$$Z:=\arg\max_{z\in\R}\left\{\sqrt{2}B(z)+A_2(z)-z^2\right\}\,,$$
and $B$ is a standard Brownian motion and $A_2$ is an independent Airy$_2$ process. The maximizer $Z$ describes the limit fluctuations of semi-infinite geodesics \cite{Pi}. To complete the proof one needs to show that $t^{-2/3}X_t$ and $t^{-2/3}Y_t$ will localize around the origin with high probability. 
 
\subsection{Local Fluctuations of an Airy Sheet}
Define 
$$H_n(u,v) = \frac{L_{2^{2/3} un^{2/3}}[2^{2/3} vn^{2/3}]_n-4 n+2^{4/3}(v-u)^2n^{1/3}}{2^{4/3} n^{1/3}}\,,(u,v)\in[0,C]^2\,.$$
The coupling method is suitable to prove relative compactness of $\{H_n(\cdot,\cdot)\,:\,n\geq 1\}$ as a collection of random two-dimensional cadlag scalar fields. It is believed that there exists a unique limiting object, called Airy sheet, and that this object is at the center of variational formulas involving the KPZ fixed point \cite{CoQuRe,MaQuRe}. Although the local comparison method does not provide uniqueness, it can be used to prove that any weak limit point is locally an additive Brownian motion with diffusion coeficient $2$.

\begin{thm}\label{mainSheet}
The collection $\{H_n(\cdot,\cdot)\,:\,n\geq 1\}$ is relatively compact in the Skorohod topology of cadlag functions on compact sets,  and any weak limit point is continuous almost surely. For a weak limit point $A(\cdot,\cdot)$ of $\{H_n(\cdot,\cdot)\,:\,n\geq 1\}$ define 
$$A^{\epsilon}(x,y):=\epsilon^{-1/2}\left(A(\epsilon x,\epsilon y)-A(0,0)\right)\,.$$ 
Then       
\begin{equation}\label{loconv}
\lim_{\epsilon\downarrow 0}A^{\epsilon}(\cdot,\cdot)\stackrel{dist.}{=} \sqrt{2}\left(B_1(\cdot)+B_2(\cdot)\right)\,,
\end{equation}
in the topology of continuous functions on compact sets, where $(B_1(x)\,,x\geq 0)$ and $(B_2(y)\,,y\geq 0)$ are two independent standard Brownian motions.
\end{thm}

\begin{rem}\label{subKPZSh}
A two parameters sub-KPZ local fluctuation sheet is defined by setting
$$\Delta_{\gamma,n}(u,v):=\frac{L_{un^{\gamma}}[vn^{\gamma}]_n-L_0[0]_n}{2^{3/2}n^{\gamma/2}}\,,u\in[0,C]\,,$$ 
where fixing $\gamma\in(0,2/3)$. The script to prove Theorem \ref{mainSheet} can be adapted to show that,
\begin{equation}\label{subloconv}
\lim_{n\to\infty}\Delta_{\gamma,n}(\cdot,\cdot)\stackrel{dist.}{=} B_1(\cdot)+B_2(\cdot)\,,
\end{equation}
in the topology of in the Skorohod topology of cadlag functions on compact sets.
\end{rem}

\section{Exit Points and Local Comparison} 
Given two profiles $\rb_1$ and $\rb_2$, the basic coupling $(L^{\rb_1},L^{\rb_2})$ is constructed by setting
$$L^{\rb_i}(\xx):=\max_{k\in C_{\xx}} \left\{ \rb_i(k) + L_k(\xx)\right\}\,.$$ 
(Recall that $L$ is a function of the same environment $\omega$). In \cite{CaPi}, the key result to study the local fluctuations of $L$ is given as follows (see Lemma 1 in \cite{CaPi} for the LPP-Poissonian version). Let $k\leq l$ and $n\geq 1$. If $Z^{\rb}[k]_n\geq0$ then
$$L[l]_n-L[k]_n\leq L^{\rb}[l]_n-L^{\rb}[k]_n\,,$$
and if $Z^{\rb}[l]_n\leq 0$ then 
$$L[l]_n-L[k]_n\geq L^{\rb}[l]_n-L^{\rb}[k]_n\,.$$
As we noted before, for the wedge profile we have that $L^\w=L$ and $Z^\w(\xx)\equiv 0$, and the above inequalities can be seen as a comparison of local increments of the LPP model with respect to $\w$ and $\rb$, in terms of the relative locations of the respective exit-points. This local comparison property can be generalize for the coupling $(L^{\rb_1},L^{\rb_2})$, with arbitrary $\rb_1$ and $\rb_2$, as follows. 
\begin{lem}\label{comparison}
Let $k\leq l$ and $n\geq 1$. If $Z^{\rb_1}[l]_n\leq Z^{\rb_2}[k]_n$ then
$$L^{\rb_1}[l]_n-L^{\rb_1}[k]_n \leq L^{\rb_2}[l]_n-L^{\rb_2}[k]_n\,.$$
\end{lem}

\noindent{\bf Proof\,\,} 
Let $\pi(\xx,\yy)$, for $\xx\leq \yy$, denote the path which attains the last-passage percolation time:
$$\pi(\xx,\yy)=\arg\max_{\pi\in\Pi(\xx,\yy)}\omega(\pi)\,,$$
so that $L(\xx,\yy)=\omega(\pi(\xx,\yy))$. Notice that
$$L(\xx,\yy)=L(\xx,\zz)+L(\zz,\yy)\,,$$
for any $\zz\in\pi(\xx,\yy)$. Denote $z_1\equiv Z^{\rb_1}[l]_n$ and $z_2\equiv Z^{\rb_2}[k]_n$ so that, by assumption, $z_1\leq z_2$. Let $\cc$ be a crossing between $\pi([z_1]_0,[l]_n)$ and $\pi([z_2]_0,[k]_n)$. Such a crossing always exists because $k\leq l$ and $z_1\leq z_2$. We remark that, by superaddivity,
$$L^{\rb_2}[l]_n  \geq  \rb_2(z_2) + L([z_2]_0,[l]_n) \geq  \rb_2(z_2) + L([z_2]_0,\cc) + L(\cc,[l]_n)\,.$$
We use this, and that (since $\cc\in\pi([z_2]_0,[k]_n)$)
$$ \rb_2(z_2) + L([z_2]_0,\cc)-L^{\rb_2} [k]_n= -L(\cc,[k]_n)\,,$$
in the following inequality:
\begin{eqnarray*}
L^{\rb_2}[l]_n - L^{\rb_2}[k]_n & \geq & \rb_2(z_2)+L([z_2]_0,\cc) + L(\cc,[l]_n) - L^{\rb_2}[k]_n\\
& = & L(\cc,[l]_n) - L(\cc,[k]_n)\,.
\end{eqnarray*}
By superaddivity,
$$ - L(\cc\,,\,[k]_n)\geq L^{\rb_1}(\cc)-L^{\rb_1}[k]_n\,,$$
and hence (since $\cc\in\pi([z_1]_0,[l]_n)$)
\begin{eqnarray*}
L^{\rb_2}[l]_n - L^{\rb_2}[k]_n & \geq & L(\cc,[l]_n) - L(\cc,[k]_n)\\
& \geq & L(\cc,[l]_n) + L^{\rb_1}(\cc)-L^{\rb_1}[k]_n\\
& = & L^{\rb_1}[l]_n-L^{\rb_1}[k]_n\,.
\end{eqnarray*}

\hfill$\Box$\\ 
 
Now the aim is to compare the local increment corresponding to a given profile $\rb$ with the local increment of the stationary profile $\s_\rho$. From now on we will put a superscript $\rho$ for quantities related to the stationary measure, such as $L^\rho\equiv L^{\s_{\rho}}$ and  $Z^\rho\equiv Z^{\s_{\rho}}$. Due to the (space) stationarity of the increments of $\s_\rho$, the location of the exit-point satisfies:
$$Z^\rho[k+h]_n \stackrel{dist.}{=} Z^\rho[k]_n +h\,.$$
To control the fluctuations of $Z^\rho$ one has to look at the so called characteristic of the system given by the direction $(a_\rho,1)$, where
$$a_\rho = \left(\frac{1-\rho}{\rho}\right)^2\,.$$
On the anti-diagonal $x+y=2n$, this corresponds to $x=n+k$ where
$$k=b_\rho n=\left( \frac{a_\rho-1}{a_\rho+1}\right) n\,.$$ 
Along this direction the exit-point oscillates around the origin in the scale   $n^{2/3}$. 
\begin{lem}\label{stcontrol} 
There exists a constant $c_1>0$ such that, for all $\rho\in [1/4,3/4]$ 
$$\limsup_{n\to\infty}\P\left(|Z^\rho[b_\rho n]_n|\geq rn^{2/3} \right)\leq \frac{c_1}{r^3}\,,$$
for all $r>1$.  
\end{lem}

\noindent{\bf Proof\,\,} If we denote by $\tilde Z^\rho$ the exit-point with respect to the positive coordinate axis then 
$$\left\{ |Z^\rho[b_\rho n]_n|\geq rn^{2/3} \right\} \subseteq \left\{ |\tilde Z^\rho[b_\rho n]_n|\geq rn^{2/3} \right\}\,.$$
Thus, Lemma \ref{stcontrol} follows from Theorem 2.2 in \cite{BaCaSe}.

\hfill$\Box$\\ 

The reason we need uniform control for all $\rho\in[1/4,3/4]$ is that we will scale $\rho$ with $n$ in a neighborhood of $1/2$. For this neighborhood we have $0\leq a_\rho\leq 9$ and so
$$0\leq\frac{1}{10}(a_\rho-1)\leq b_\rho \,,\,\mbox{ for }\,\rho\in [1/4,1/2]\,,$$  
and 
$$0\geq\frac{1}{10}(a_\rho-1)\geq b_\rho\,,\,\mbox{ for }\,\rho\in [1/2,3/4]\,.$$  
We note that $a_\rho$ is a decreasing function of $\rho$ and its derivative is bounded in the interval $[1/4,3/4]$. So by using the mean value theorem, one can see that  there are constants $c_2,c_3\in(0,\infty)$ such that 
$$0\leq -c_2\left(\rho-\frac{1}{2}\right)\leq b_\rho \,,\,\mbox{ for }\,\rho\in [1/4,1/2]\,,$$  
and 
$$0\geq-c_3\left(\rho-\frac{1}{2}\right)\geq b_\rho\,,\,\mbox{ for }\,\rho\in [1/2,3/4]\,.$$
We will be interested in the regime 
$$\rho_n^{\pm}=1/2\pm\frac{r}{n^{1/3}}\,,$$ 
where $r$ is fixed and $n$ is large enough so that $\rho^\pm_n\in[1/4,3/4]$. Therefore, 
$$0\leq c_2rn^{2/3}\leq b_{\rho_n^-}n \,,$$  
and 
$$0\geq -c_3rn^{2/3}\geq b_{\rho_n^+}n\,.$$

\begin{lem}\label{compcontrol}
For $r,C\geq 0$ set 
$$\rho_n^{\pm}=\rho_n^{\pm}(r)=1/2\pm\frac{r}{n^{1/3}}\,,$$
and define the event
$$E_n(r):=\left\{Z^{\rho_n^-}[Cn^{2/3}]_n\leq Z^\rb [0]_n\,\,\mbox{ and }Z^{\rb}[Cn^{2/3}]_n\leq Z^{\rho_n^+} [0]_n\,\, \right\}\,.$$
Under Assumption \ref{exitcontrol},  
$$\limsup_{r\to\infty}\limsup_{n\to\infty}\P\left(E_n(r)^c\right)=0\,.$$
\end{lem}

\noindent{\bf Proof\,\,} 
Let us first show that 
$$\limsup_{r\to\infty}\limsup_{n\to\infty}\P\left( Z^{\rb}[Cn^{2/3}]_n> Z^{\rho_n^+} [0]_n\right)=0\,. $$ 
Since,  
$$\P\left( Z^{\rb}[Cn^{2/3}]_n> Z^{\rho_n^+} [0]_n\right)\leq\P\left( |Z^{\rb}[Cn^{2/3}]_n|>\frac{c_3}{2}r n^{2/3}\right)+\P\left(Z^{\rho_n^+} [0]_n<\frac{c_3}{2}r n^{2/3}\right)\,,$$
By Assumption \ref{exitcontrol}, we only need to control 
$$ \limsup_{n\to\infty}\P\left(Z^{\rho_n^+} [0]_n<\frac{c_3}{2}r n^{2/3}\right)\,.$$
On the other hand, $Z^{\rho_n^+} [0]_n \stackrel{dist.}{=}Z^{\rho_n^+} [b_{\rho_n^+}n]_n-b_{\rho_n^+}n$, and so 
\begin{eqnarray*}
\P\left(Z^{\rho_n^+} [0]_n<\frac{c_3}{2}r n^{2/3}\right)&=&\P\left(Z^{\rho_n^+} [b_{\rho_n^+}n]_n<\frac{c_3}{2}r n^{2/3}+b_{\rho_n^+}n\right)\\
 &\leq &\P\left(Z^{\rho_n^+} [b_{\rho_n^+}n]_n<\left(\frac{c_3}{2}-c_3\right)rn^{2/3}\right)\\
 &=& \P\left(Z^{\rho_n^+} [b_{\rho_n^+}n]_n<-\frac{c_3}{2}rn^{2/3}\right)\\
  &\leq& \P\left(|Z^{\rho_n^+} [b_{\rho_n^+}n]_n|>\frac{c_3}{2}rn^{2/3}\right)\,.
\end{eqnarray*}
By Lemma \ref{stcontrol}, this shows that
\begin{equation}\label{+control}
\limsup_{n\to\infty}\P\left(Z^{\rho_n^+} [0]_n<\frac{c_3}{2} n^{2/3}\right)\leq \frac{c'_1}{r^{3}}\,,
\end{equation}
for large enough $r$. To estimate  
$$\limsup_{n\to\infty}\P\left( Z^{\rb}[0]_n< Z^{\rho_n^-} [Cn^{2/3}]_n\right)\,,$$
we use a similar argument:
$$\P\left( Z^{\rb}[0]_n< Z^{\rho_n^-} [Cn^{2/3}]_n]_n\right)\leq \P\left( |Z^{\rb}[0]_n|>\frac{c_2}{2}r n^{2/3}\right)+ \P\left(Z^{\rho_n^-} [Cn^{2/3}]_n]>-\frac{c_2}{2}r n^{2/3}\right)\,,$$ 
and
\begin{eqnarray*}
\P\left(Z^{\rho_n^-} [Cn^{2/3}]_n>-\frac{c_2}{2}r n^{2/3}\right)&=&\P\left(Z^{\rho_n^-} [b_{\rho_n^-}n]_n>-\frac{c_2}{2}r n^{2/3}+b_{\rho_n^-}n - Cn^{2/3}\right)\\
 &\leq &\P\left(Z^{\rho_n^-} [b_{\rho_n^-}n]_n>\left(\frac{c_2}{2}r-C\right)n^{2/3}\right)\,.
\end{eqnarray*}
Hence (by Lemma \ref{stcontrol}), for large enough $r$,
\begin{equation}\label{-control}
\limsup_{n\to\infty}\P\left(Z^{\rho_n^-} [Cn^{2/3}]_n>-\frac{c_2}{2}r n^{2/3}\right)\leq \frac{c''_1}{r^{3}}\,.
\end{equation}

\hfill$\Box$\\ 

The stationary profile with parameter $\rho$ has mean
$$\mu_\rho=\E\left(X(\rho)\right)=\E\left(\Exp_1(1-\rho)\right)-\E\left(\Exp_2(\rho)\right)=\frac{2\rho-1}{\rho(1-\rho)}\,,$$  
and 
$$\mu_{\rho_n^\pm}n^{1/3}=\pm\frac{r}{\rho_n^{\pm}(1-\rho_n^{\pm})}\,.$$
Hence (minimize the denominator),  
\begin{equation}\label{compest1}
\mu_{\rho^+_n}n^{1/3}\leq 6r\,\,\mbox{ and }\,\,\mu_{\rho^-_n}n^{1/3}\geq -6r\,,
\end{equation}
for fixed $r>0$, and large enough $n$ such that $\rho_n^\pm\in[1/4,3/4]$. Let 
$$B^{\pm}_n(u):=\frac{L^{\rho_n^\pm}[un^{2/3}]_n-L^{\rho_n^\pm}[0]_n-\mu_{\rho_n^\pm}un^{2/3}}{2^{3/2}n^{1/3}}\,.$$

\begin{lem}\label{compest}
On the event $E_n(r)$,   
$$B^{-}_n(v)-B^{-}_n(u)-3\sqrt{2}(v-u)r \leq \Delta^\rb_n(v)-\Delta^\rb_n(u)\leq B^{+}_n(v)- B^{+}_n(u)+3\sqrt{2}(v-u)r\,,$$
for all $u,v\in[0,C]$.
\end{lem}

\noindent{\bf Proof\,\,} For fixed $n\geq 1$, $Z^\rb[k]_n$ is a nondecreasing function of $k$ and, on the event $E_n(r)$,
$$Z^{\rho_n^-}[vn^{2/3}]_n\leq Z^{\rho_n^-}[Cn^{2/3}]_n\leq Z^\rb [0]_n\leq Z^\rb [un^{2/3}]_n\,,$$
and 
$$Z^{\rb}[vn^{2/3}]_n\leq Z^{\rb}[Cn^{2/3}]_n\leq Z^{\rho_n^+} [0]_n \leq Z^{\rho_n^+} [un^{2/3}]_n\,,$$
for all $u,v\in[0,C]$. By Lemma \ref{comparison}, on the event $E_n(r)$,
$$L^{\rho_n^-}[vn^{2/3}]_n-L^{\rho_n^-}[un^{2/3}]_n\leq L^\rb[vn^{2/3}]_n-L^\rb[un^{2/3}]_n\leq L^{\rho_n^+}[vn^{2/3}]_n-L^{\rho_n^+}[un^{2/3}]_n\,,$$
which yields
$$B^{-}_n(v)-B^{-}_n(u)+2^{-3/2}\mu_{\rho^-_n}n^{1/3}(v-u)\leq \Delta_n^\rb(v)-\Delta_n^\rb(u)\,,$$
and 
$$\Delta_n^\rb(v)-\Delta_n^\rb(u)\leq B^{+}_n(v)-B^{+}_n(u)+2^{-3/2}\mu_{\rho^+_n}n^{1/3}(v-u)\,,$$
for all $u,v\in[0,C]$. Together with \eqref{compest1}, this proves Lemma \ref{compest}

\hfill$\Box$\\ 

The modulus of continuity of a process $(X(u)\,,\,u\in [0,C]$ is defined as 
$$W_X(\delta):=\sup_{u,v\in [0,C],\,|u-v|\leq \delta} |X(u)-X(v)|\,,\mbox{ where }\delta\in[0,C]\,.$$
\begin{coro} \label{Corcompest}
On the event $E_n(r)$,
$$W_{\Delta^b_n}(\delta)\leq \max\left\{W_{B^{-}_n}(\delta)\,,\,W_{B^{+}_n}(\delta) \right\} + 3\sqrt{2}\delta r\,.$$
\end{coro}

\noindent{\bf Proof of Theorem \ref{main}.\,}
Denote
$$R_1(r):=\limsup_{n\to\infty}\P\left(E_n(r)^c\right)\,.$$
Since $\cal D$ is separable and complete (under a suitable metric) a family of probability measures on $\cal D$ is relatively compact if and only if it is tight (Prohorov's Theorem). Thus, we will use Theorem 15.5 in \cite{Bi}, which states that if: 
\begin{itemize}
\item $\limsup_{a\to\infty}\limsup_{n\to\infty}\P\left(|X_n(0)|>a\right)=0$;
\item for every $\eta>0$, $\limsup_{\delta\downarrow 0}\limsup_{n\to\infty}\P\left(W_{X_n}(\delta)>\eta\right)=0$;
\end{itemize}
then $\left\{X_{n}\,:\,n\geq 1\right\}$ is tight in the Skorohod topology of cadlag functions on compact sets, and any weak limit point is continuous almost surely. Since $\Delta^b_n(0)=0$ we only need to check the second condition: by Corollary \ref{Corcompest},
$$\P\left(W_{\Delta^b_n}(\delta)>\eta\right)\leq \P\left(W_{B^{-}_n}(\delta)>\eta-3\sqrt{2}\delta r\right)+\P\left(W_{B^{+}_n}(\delta)>\eta-3\sqrt{2}\delta r\right)+\P\left(E_n(r)^c\right)\,.$$
Clearly $\lim_{n\to\infty}B^{\pm}_n \stackrel{dist.}{=} B$, where $B$ is a standard Brownian motion. Thus,
$$\limsup_{n\to\infty}\P\left(W_{\Delta^b_n}(\delta)>\eta\right)\leq 2\P\left(W_{B}(\delta)>\eta-3\sqrt{2}\delta r\right)+R_1(r)\,.$$
If we choose $r=r_\delta:=\delta^{-\beta}$, for a fixed $\beta\in(0,1)$, then 
\begin{equation*} 
\limsup_{\delta\downarrow 0}\limsup_{n\to\infty}\P\left(W_{\Delta^b_n}(\delta)>\eta\right)\leq 2\limsup_{\delta\downarrow 0}\P\left(W_{B}(\delta)>\eta-3\sqrt{2}\delta^{1-\beta}\right)+\limsup_{\delta\downarrow 0}R_1(r_\delta)=0\,.
\end{equation*}
We use Lemma \ref{compcontrol} to control the probability of $E_n(r_\delta)$, while for the other $\limsup_{\delta\downarrow 0}$ we just use use that if $X$ is a stochastic process in the space of continuous functions (uniform topology) then  
$$\limsup_{\delta \downarrow 0}\P\left(W_{X}(\delta)>\eta\right)=0\,.$$

To prove \eqref{loconv}, given a process $X$, we denote
$$ X^\epsilon\,:\, u\mapsto\epsilon^{1/2}X(\epsilon^{-1}u)\,.$$ 
By Lemma \ref{compest}, if $|v-u|\leq \delta$ then,
\begin{equation}\label{localBM}
B^{-,\epsilon}_n(v)-B^{-,\epsilon}_n(u)-3\sqrt{2}\epsilon^{1/2}\delta r \leq \Delta^{\rb,\epsilon}_n(v)-\Delta^{\rb,\epsilon}_n(u)\leq B^{+,\epsilon}_n(v)- B^{+,\epsilon}_n(u)+3\sqrt{2}\epsilon^{1/2}\delta r\,.
\end{equation}
Since $\epsilon^{1/2}B(\epsilon^{-1}x) \stackrel{dist.}{=} B(x)$, if $\Delta^\rb$ is any weak limit point of $\left\{\Delta_n^\rb\,:\,n\geq 1\right\}$, then    
$$\P\left(W_{\Delta^{\rb,\epsilon}}(\delta)>\eta\right)\leq 2\P\left(W_{B}(\delta)>\eta-3\sqrt{2}\delta\epsilon^{1/2} r\right)+R_1(r)\,.$$
If we now set $r_\epsilon:=\epsilon^{-\beta}$, with $\beta\in(0,1/2)$, then 
$$\limsup_{\epsilon\downarrow 0}\P\left(W_{\Delta^{\rb,\epsilon}}(\delta)>\eta\right)\leq 2\P\left(W_{B}(\delta)>\eta\right)\,,$$
and hence
$$\limsup_{\delta\downarrow 0}\limsup_{\epsilon\downarrow 0}\P\left(W_{\Delta^{\rb,\epsilon}}(\delta)>\eta\right)= 0\,,$$
which implies tightness of $\Delta^{\rb,\epsilon}$ (in the space of continuous processes). By \eqref{localBM}, if $u_i\in[0,A]$ then 
$$\P\left(\cap_{i=1}^j\left\{ \Delta^{\rb,\epsilon}(u_i)\leq x_i\right\}\right)\leq \P\left(\cap_{i=1}^j\left\{ B(u_i)\leq x_i +3\sqrt{2}A\epsilon^{1/2-\beta}\right\}\right)+R_1(r_\epsilon)\,  $$
and
$$\P\left(\cap_{i=1}^j\left\{ \Delta^{\rb,\epsilon}(u_i)\leq x_i\right\}\right)\geq \P\left(\cap_{i=1}^j\left\{ B(u_i)\leq x_i -3\sqrt{2}A\epsilon^{1/2-\beta}\right\}\right)+R_1(r_\epsilon)\,,$$
which proves finite-dimensional convergence:
$$\lim_{\epsilon\downarrow 0}\P\left(\cap_{i=1}^j\left\{ \Delta^{\rb,\epsilon}(u_i)\leq x_i\right\}\right)= \P\left(\cap_{i=1}^j\left\{ B(u_i)\leq x_i \right\}\right)\,.$$

To prove \eqref{subloconv} one needs to pick $\gamma'\in(\gamma,2/3)$, set
$$\rho_n^{\pm}=1/2\pm\frac{1}{n^{\gamma'/2}}\,,$$ 
and define the exit-point comparison event  
$$E_n:=\left\{Z^{\rho_n^-}[Cn^{\gamma}]_n\leq Z^\rb [0]_n\,\,\mbox{ and }Z^{\rb}[Cn^\gamma]_n\leq Z^{\rho_n^+} [0]_n\,\, \right\}\,.$$
By using $r_n=1/3-\gamma'/2$ as in the proof of Lemma \ref{compcontrol}, one gets that 
$$\limsup_{n\to\infty}\P\left(E_n^c\right)=0\,,$$
and the rest of the proof follows by applying local comparison mutatis mutandis.

\section{Checking Assumption 1}\label{check}
We will prove Theorem \ref{local} separately for each case. The starting point of the proof will always be the same: let $\rb$ be an initial profile and consider the event 
\begin{equation*}
\left\{Z^\rb[0]_n\geq u \right\}=\left\{\exists\,z\in[u,n]\,:\,\rb(z)+L_z[0]_n=L^{\rb}[0]_n\right\}\,.
\end{equation*}
Pick an invariant profile $\s_\rho$ and construct $L^{\s_\rho}$ and $L^{\rb}$ simultaneously using the basic coupling. Since  
$$L_z[0]_n\leq L^{\s_\rho}[0]_n - s_\rho(z)\,\mbox{ and }\,L[0]_n=L_0[0]_n+\rb(0)\leq L^\rb[0]_n\,$$
(recall that $\rb(0)=0$ for the second inequality), we have that 
\begin{equation}\label{exit1}
\left\{Z^\rb[0]_n\geq u \right\}\subseteq\left\{\exists\,z\in[u,n]\,:\,\s_\rho(z)-\rb(z)\leq L^{\s_\rho}[0]_n-L[0]_n\right\}\,.
\end{equation}

\noindent{\bf Flat Profile.\,} For the flat profile $\rb=\f\equiv 0$, we can rewrite the event in the right hand side of \eqref{exit1} as, 
$$\ra(u):=\m(u)+\s_\rho(u)\leq L^{\s_\rho}[0]_n-L[0]_n\,,$$  
where $\m(u)$ is the minimum of $\s_\rho(z)-\s_\rho(u)$ for $z\geq u$. Hence,
\begin{equation}\label{exit2}
\P\left(Z^\f[0]_n\geq u\right)\leq \P\big(\ra(u)\leq L^{\s_\rho}[0]_n-L[0]_n\big)\,. 
\end{equation}
A straightforward computation shows that, 
$$\E \ra(u)= \E \m(u)+\frac{2\rho - 1}{\rho(1-\rho)}u\,,$$
and
$$\Var \ra(u) = \Var \m(u) +\left(\frac{1}{(1-\rho)^2}+\frac{1}{\rho^2}\right)u\,$$
($\m(u)$ and $\s_\rho(u)$ are independent). The random walk $\left(\s(u)-\s(u+n)\,,\,n\geq 0\right)$ has a negative drift, and its maximum has a well known distribution \cite{Res}. This allow us to get that 
$$\E \m(u)=-\frac{1}{2\rho-1}\frac{1-\rho}{\rho}\,,$$ 
and  
$$\Var \m(u)=\frac{2(1-\rho)}{\rho(2\rho-1)^2}-\frac{(1-\rho)^2}{\rho^2(2\rho-1)^2}=\frac{1}{(2\rho-1)^2}\left(2\frac{(1-\rho)}{\rho}-\frac{(1-\rho)^2}{\rho^2}\right)\,.$$

Now, for $c>0$ (we will set its value later),
\begin{eqnarray}
\nonumber\P\big(\ra(u)\leq L^{\s_\rho}[0]_n-L[0]_n\big)&\leq & \P\left(L[0]_n-4n\leq -\frac{r^2n^{1/3}}{c}\right)\\
\label{exit3}&+&\P\left(\ra(u)\leq L^{\s_\rho}[0]_n-4n+\frac{r^2n^{1/3}}{c}\right)\,.
\end{eqnarray}
By Theorem 2.4 in \cite{BaCaSe}, for $\alpha\in(0,1)$ there exists $c_1=c_1(\alpha)>0$ such that for all $n\geq 1$ and $r>0$,
\begin{equation*}
\P\left(L[0]_n-4n\leq -\frac{r^2n^{1/3}}{c}\right)\leq \frac{c_1 c^{3\alpha/2}}{(r)^{3\alpha}}\,.
\end{equation*}
and so,   
\begin{equation}\label{wedgetail}
\limsup_{r\to\infty}\limsup_{n\to\infty}\P\left(L[0]_n-4n\leq -\frac{r^2n^{1/3}}{c}\right)=0\,.
\end{equation}

To estimate the second term in the rhs of \eqref{exit3} we do as fllows:
\begin{eqnarray}
\nonumber\P\left(\ra(u)\leq L^{\s_\rho}[0]_n-4n+\frac{r^2n^{1/3}}{c}\right)&\leq &\P\left(L^{\s_\rho}[0]_n-4n\geq \E \ra(u) - \frac{2r^2n^{1/3}}{c}\right)\\
\label{exit4}&+& \P\left(\ra(u)\leq \E \ra(u) - \frac{r^2n^{1/3}}{c}\right)\,.
\end{eqnarray}
Write
$$\P\left(L^\s[0]_n-4n\geq \E \ra(u) - \frac{2r^2n^{1/3}}{L}\right)=\P\left(L^\s[0]_n-\E L^\s[0]_n\geq \Lambda - \frac{2r^2n^{1/3}}{c}\right)\,,$$
where
\begin{eqnarray*}
\Lambda &=&4n-\E L^\s[0]_n+\E \ra(u)\\
&=& \left(4-\frac{1}{\rho(1-\rho)}\right)n+\frac{2\rho - 1}{\rho(1-\rho)}u +\E \m(u) \\
&=& \frac{(4\rho(1-\rho)-1)n+(2\rho - 1)u}{\rho(1-\rho)} +\E \m(u)\\
&\geq &\frac{(4\rho(1-\rho)-1)n+(2\rho - 1)u}{4}+\E \m(u)\,.
\end{eqnarray*}
By maximizing the term in the numerator, we get 
$$\rho(u,n)=\frac{1}{2}+\frac{u}{4n}\,,$$
and we set $u:=rn^{2/3}$. For this choice of $\rho(u,n)$, we have
$$ \frac{(4\rho(1-\rho)-1)n+(2\rho - 1)u}{4}=\frac{u^2}{16 n}=\frac{r^2}{16}n^{1/3}\,.$$
We still have to control $\E \m(u)$:
$$\E \m(u)=-\frac{1}{2\rho-1}\frac{1-\rho}{\rho}\geq-\frac{4}{r}n^{1/3}\,.$$ 
Thus
$$\Lambda \geq \left(\frac{r^2}{16}-\frac{4}{r}\right)n^{1/3}\geq \frac{r^2}{32}n^{1/3}\,,$$
for $r>128$. If we choose $c=128$, then we finally get
\begin{eqnarray}
\nonumber\P\left(L^\s[0]_n-4n\geq \E \ra(u) - \frac{2r^2n^{1/3}}{128}\right)&=&\P\left(L^\s[0]_n-\E L^\s[0]_n\geq \Lambda - \frac{r^2n^{1/3}}{64}\right)\\
\nonumber &\leq & \P\left(L^\s[0]_n-\E L^\s[0]_n\geq \frac{r^2n^{1/3}}{64}\right)\\ 
\nonumber &\leq &(64)^2\frac{\Var L^{\s}[0]_n}{r^4n^{2/3}}\,.
\end{eqnarray}
By Lemma 4.7 in \cite{BaCaSe},
$$\Var L^{\s_\rho}[0]_n\leq \frac{ \Var L^{\s_{1/2}}[0]_n }{4\rho^2}+n\frac{2\rho -1}{\rho^2(1-\rho)^2}\,.$$
Together with Theorem 2.1 in \cite{BaCaSe}, this shows that 
$$\limsup_{n\to\infty}\frac{\Var L^{\s_\rho}[0]_n}{n^{2/3}}\leq\limsup_{n\to\infty}\frac{\Var L^{\s_{1/2}}[0]_n}{n^{2/3}} + 4r\leq c_2+ 4r\,,$$  
for some constant $c_2>0$. Therefore,
\begin{equation}\label{exit5}
\limsup_{n\to\infty}\P\left(L^\s[0]_n-4n\geq \E \ra(u) - \frac{2r^2n^{1/3}}{128}\right)\leq (64)^2\left(\frac{c_2}{r^4}+\frac{4}{r^3}\right)\,.
\end{equation}

Now, for $u=rn^{2/3}$,
$$\lim_{n\to\infty}\Var X(\rho)= 8\,\mbox{ and }\,\lim_{n\to\infty}\frac{\Var \m(u)}{n^{2/3}}= \frac{1}{16r^{2}}\,,$$
which implies that  
\begin{eqnarray}
\nonumber\limsup_{n\to\infty}\P\left(\ra(u)\leq \E \ra(u) - \frac{r^2n^{1/3}}{128}\right)&\leq & \limsup_{n\to\infty}\frac{(128)^2}{r^4n^{2/3}}\Var \ra(u)\\
\nonumber&= &\limsup_{n\to\infty}\frac{(128)^2}{r^4n^{2/3}}\left(\Var \m(u) + (\Var X) u \right)\\
\label{exit6}&=& \frac{(128)^2}{16r^{6}}+\frac{8(128)^2}{r^3}\,.
\end{eqnarray}

From \eqref{exit2}, \eqref{exit3},\eqref{wedgetail}, \eqref{exit4}, \eqref{exit5} and \eqref{exit6}, we see that  
\begin{equation}\label{exit7}
\limsup_{r\to\infty}\limsup_{n\to\infty}\P\left(Z^\f[0]^+_n\geq rn^{2/3}\right)=0\,.
\end{equation}
By translation invariance, $Z^\f[h]_n\stackrel{dist.}{=}Z^{\f}[0]+h$, and by symmetry $Z^{\f}[0]^-_n\stackrel{dist.}{=}Z^{\f}[0]^+_n$, and so \eqref{exit7} implies Assumption \ref{exitcontrol}.

\hfill$\Box$\\

\noindent\paragraph{\bf Mixed Profiles.} The proof of Assumption \ref{exitcontrol} for mixed profiles is based on the same ideas as befored. We just need to adapt a few details in each case. 

\noindent\paragraph{\bf Wedge to Flat.} In that case $Z^{\w\f}[k]_n\geq 0$. By \eqref{exit1}, 
\begin{equation*}
\P\left(Z^{\w\f}[\pm n^{2/3}]_n\geq u\right)\leq \P\big(\ra(u)\leq L^{\s_\rho}[\pm n^{2/3}]_n-L[\pm n^{2/3}]_n\big)\,, 
\end{equation*}
with $\ra(u)$ as define before. By writing 
$$L[\pm n^{2/3}]_n-4n=\left(L[0]_n-4n\right)+\left(L[\pm n^{2/3}]_n-L[0]_n\right)\,,$$
we can use Theorem 2.4 in \cite{BaCaSe} to deal with the first term, and tightness of the local increment \cite{CaPi} to deal with the second one, and get that   
\begin{equation}\label{tailTW}
\limsup_{r\to\infty}\limsup_{n\to\infty}\P\left(L[\pm n^{2/3}]_n-4n\leq -\frac{r^2n^{1/3}}{128}\right)=0\,.
\end{equation}

Now, for the same choice of $\rho(u,n)$ as before, this leaves us with (recall the proof of \eqref{exit4}, \eqref{exit5} and \eqref{exit6})
\begin{eqnarray}
\nonumber\P\left(\ra(u)\leq L^{\s_\rho}[\pm n^{2/3}]_n-4n+\frac{r^2n^{1/3}}{128}\right)&\leq &\P\left(L^{\s_\rho}[\pm n^{2/3}]_n-\E L^{\s_\rho}[0]_n\geq \frac{r^2}{64}n^{1/3}\right)\\
\nonumber &+& \P\left(\ra(u)\leq \E \ra(u) - \frac{r^2n^{1/3}}{128}\right)\\
\nonumber &= &\P\left(L^{\s_\rho}[0]_n-\E L^{\s_\rho}[0]_n+\s_\rho(\pm n^{2/3})\geq\frac{r^2}{64}n^{1/3} \right)\\
\nonumber &+& \P\left(\ra(u)\leq \E \ra(u) - \frac{r^2n^{1/3}}{128}\right)\\
\nonumber &= &\P\left(L^{\s_\rho}[0]_n-\E L^{\s_\rho}[0]_n\geq \frac{r^2}{128}n^{1/3}\right)\\
\nonumber &+& \P\left( \s_\rho(\pm n^{2/3})\geq \frac{r^2}{128}n^{1/3}\right)\\
\nonumber &+& \P\left(\ra(u)\leq \E \ra(u) - \frac{r^2n^{1/3}}{128}\right)\,,
\end{eqnarray}
where we use that $L^{\s_\rho}[\pm n^{2/3}]_n\stackrel{dist.}{=}L^{\s_\rho}[0]_n+\s_\rho(\pm n^{2/3})$. The only term that did not appear before is 
\begin{equation*}
\P\left( \s_\rho(\pm n^{2/3})\geq \frac{r^2}{128}n^{1/3}\right)=\P\left( \s_\rho(\pm n^{2/3})-\E\s_\rho(\pm n^{2/3}) \geq\frac{r^2}{128}n^{1/3}-\E\s_\rho(\pm n^{2/3})\right)\,. 
\end{equation*}
But since $\s_{\rho}$ is a random walk with drift 
$$\E\s_\rho(\pm n^{2/3})=\frac{2\rho-1}{\rho(1-\rho)}n^{2/3}\sim r n^{1/3}\,,$$
it is not hard to see that 
$$\limsup_{r\to\infty}\limsup_{n\to\infty}\P\left( \s_\rho(\pm n^{2/3})\geq \frac{r^2}{128}n^{1/3}\right)=0\,.$$

\noindent\paragraph{\bf Wedge to Stationary.} Again we have $Z^{\w\s}[k]_n\geq 0$  and, for $\rho>1/2$,
\begin{equation*}
\P\left(Z^{\w\s}[\pm n^{2/3}]_n\geq u\right)\leq \P\big(\ra(u)\leq L^{\s_\rho}[\pm n^{2/3}]_n-L[\pm n^{2/3}]_n\big)\,, 
\end{equation*}
but now $\ra(u)$ has a different form:
$$\ra(u):=\m(u)+\rt_\rho(u)\,,$$
where 
$$\rt_\rho(u):=\s_\rho(u)-\s_{1/2}(u)\,$$
and $\m(u)$ is the minimum of the random walk 
\begin{equation}\label{rw}
(\s_\rho(u+n)-\s_{1/2}(u+n))-\rt_\rho(u))\mbox{ for }n\geq 0\,.
\end{equation}
Following \cite{BaCaSe}, we couple 
$$\s_\rho(k)=\sum_{i=1}^k\bar\zeta_i\,\,\mbox{ with }\,\,\s_{1/2}(k)=\sum_{i=1}^k\zeta_i\,,$$ by taking  
$$\zeta_i:=\Exp_{1,i}(1/2)-\Exp_{2,i}(1/2)\,\,\mbox{ and }\,\,\bar\zeta_i:=\frac{1}{2(1-\rho)}\Exp_{1,i}(1/2)-\frac{1}{2\rho}\Exp_{2,i}(1/2)\,.$$ 
The random walk \eqref{rw} will have increments 
$$\frac{2\rho-1}{2(1-\rho)}\Exp_{1,i}(1/2)+\frac{2\rho-1}{2\rho}\Exp_{2,i}(1/2)\,,\mbox{ for }i\geq 1\,,$$
and thus $\m(u)=0$. Therefore, we are left with 
\begin{equation*}
\P\left(Z^{\w\s}[\pm n^{2/3}]_n\geq u\right)\leq \P\big(\rt_\rho(u)\leq L^{\s_\rho}[\pm n^{2/3}]_n-L[\pm n^{2/3}]_n\big)\,. 
\end{equation*}
We can apply \eqref{tailTW} again and get an inequality similar to \eqref{exit3}. We note that $\rt_\rho(u)$ has the same mean as $\s_\rho(u)$, so similar estimates can be done to get the analog of \eqref{exit5}. Since, $2\rho-1= u(4n)^{-1}$ and $u= rn^{2/3}$,  
$$\Var\rt_\rho(u)\leq c_3(2\rho-1)^2u\sim \frac{u^3}{n^2}\sim r^3\,,$$
we get the analog of \eqref{exit6} but with $\limsup_{n}\left(\cdot\right)=0$.   

\noindent\paragraph{\bf Flat to Stationary.} We can estimate $Z^{\f\s}[k]^-_n$ by doing the same reasoning applied to $Z^{\w\f}[k]_n$. For $k<0$, the stationary measure has increment $\Exp_{1}(\rho)-\Exp_{2}(1-\rho)$,        
and we will have to chose $\rho(u,n)=1-u(4n)^{-1}$ instead. For $Z^{\f\s}[k]^+_n$ we can proceed exactly as we did to analise $Z^{\w\s}[k]_n$.

\section{Local Comparison in the Airy Sheet Context}
The proof of Theorem \ref{mainSheet} is again based on local comparison. The proof of tightness is very similar to what we have seen so far. Define
$$\Delta_n(u,v):=\frac{L_{un^{2/3}}[vn^{2/3}]_n-L[0]_n}{2^{3/2}n^{1/3}}\,,u\in[0,C]\,.$$
Thus,
$$H_{n}(u,v)=H_{n}(0,0)+2^{1/6}\Delta_n(2^{2/3}u,2^{2/3} v)+(v-u)^2\,,$$
and we only need to analise tightness of $\Delta_n(\cdot,\cdot)\,,\,n\geq 1$, and local convergence of a weak limit point. 
To analise the modulus of continuity of a random scalar field $X:\R^2\to\R$, 
$$W_X(\delta):=\sup_{x,y\in [0,C]^2,\,|x-y|\leq \delta} |X(x)-X(y)|\,,$$
we write
$$\P\left(W_{\Delta_n}(\delta)>\eta\right)\leq\P\left(W_{\Gamma^1_n}(\delta)>\eta/2\right)+\P\left(W_{\Gamma^2_n}(\delta)>\eta/2\right)\,,$$
where
$$\Gamma^1_n(u_1,v_1,v_2):=\frac{L_{u_1n^{2/3}}[v_1n^{2/3}]_n-L_{u_1n^{2/3}}[v_2n^{2/3}]_n}{2^{3/2}n^{1/3}}\,,$$
and 
$$\Gamma^2_n(u_1,u_2,v_2):=\frac{L_{u_1n^{2/3}}[v_2n^{2/3}]_n-L_{u_2n^{2/3}}[v_2n^{2/3}]_n}{2^{3/2}n^{1/3}}\,,$$
so that 
$$\Delta_n(u_1,v_1)-\Delta_n(u_2,v_2)=\Gamma^1_n+\Gamma^2_n\,.$$
We will prove that 
\begin{equation}\label{modG1}
\limsup_{\delta\downarrow 0}\limsup_{n\to\infty}\P\left(W_{\Gamma^1_n}(\delta)>\eta\right)=0\,,
\end{equation}
and it should be clear that the same type of argument implies the analog result for $\Gamma_2$ (by interchanging the rules of $u_1$ and $v_2$). We start by proving two   lemmas that are similar to Lemma \ref{compcontrol} and Lemma \ref{compest}, respectively.
\begin{lem}\label{compcontrol1}
For $r,C\geq 0$ set 
$$\rho_n^{\pm}=\rho_n^{\pm}(r)=1/2\pm\frac{r}{n^{1/3}}\,,$$
and define the event
$$\bar E_n(r):=\left\{Z^{\rho_n^-}[Cn^{2/3}]_n\leq 0\,\,\mbox{ and }Z^{\rho_n^+} [0]_n\geq Cn^{2/3}\,\, \right\}\,.$$
Then,  
$$\limsup_{r\to\infty}\limsup_{n\to\infty}\P\left(\bar E_n(r)^c\right)=0\,.$$
\end{lem}

\noindent{\bf Proof\,\,} We follow the same argument used in Lemma \ref{compcontrol},
$$\P\left(\bar E_n(r)^c\right)\leq \P\left(Z^{\rho_n^-}[Cn^{2/3}]_n> 0 \right)+\P\left(Z^{\rho_n^+}[0]_n< Cn^{2/3} \right)\,,$$
and use \eqref{+control} and \eqref{-control} to bound the probabilities in the right-hand side of the above inequality.

\hfill$\Box$\\ 

\begin{lem}\label{compest2}
Let 
$$B^{\pm}_n(v):=\frac{L^{\rho_n^\pm}[vn^{2/3}]_n-L^{\rho_n^\pm}[0]_n-\mu_{\rho_n^\pm}vn^{2/3}}{2^{3/2}n^{1/3}}\,.$$
On the event $\bar E_n(r)$,   
$$B^{-}_n(v_2)-B^{-}_n(v_1)-3\sqrt{2}(v_2-v_1)r \leq \Gamma^1_n(u_1,v_1,v_2)\leq B^{+}_n(v_2)- B^{+}_n(v_1)+3\sqrt{2}(v_2-v_1)r\,,$$
for all $u_1\in[0,C]$ and $v_1,v_2\in[0,C]$.
\end{lem}

\noindent{\bf Proof\,\,} It follows from Lemma \ref{comparison}, as in the proof of  Lemma \ref{compest}. To see this, consider the profile 
$$\bar\w(k)=\left\{\begin{array}{ll}0 &\mbox{ for } k=u_1n^{2/3}\\
-\infty& \mbox{ for } k\neq u_1n^{2/3}\,,\end{array}\right.$$ 
and notice that, on the event $\bar E_n(r)$, $Z^{\rho_n^-}[Cn^{2/3}]_n\leq u_1n^{2/3}$ and $Z^{\rho_n^+}[0]_n\geq u_1n^{2/3}$, for all $u_1\in[0,C]$. 

\hfill$\Box$\\ 

Therefore, by Lemma \ref{compest1}, on the event $\bar{E}_n(r)$,
$$W_{\Gamma^1_n}(\delta)\leq \max\left\{W_{B^{-}_n}(\delta)\,,\,W_{B^{+}_n}(\delta) \right\} + 3\sqrt{2}\delta r\,,$$
and hence
$$\P\left(W_{\Gamma^1_n}(\delta)>\eta\right)\leq \P\left(W_{B^{-}_n}(\delta)>\eta-3\sqrt{2}\delta r\right)+\P\left(W_{B^{+}_n}(\delta)>\eta-3\sqrt{2}\delta r\right)+\P\left(\bar E_n(r)^c\right)\,.$$
Since $\lim_{n\to\infty}B^{\pm}_n \stackrel{dist.}{=} B$, where $B$ is a standard Brownian motion, together with Lemma \ref{compcontrol1} this implies \eqref{modG1}.
\newline

To prove tightness of $\epsilon^{-1/2}\Delta(\epsilon u,\epsilon v)$ one just need to repeat the same steps as in \eqref{localBM}, but now  using Lemma \ref{compcontrol1} and Lemma \ref{compest2} (and its counter part for $\Gamma^2_n$). It is not hard to guess that it will converge to a sum of two Brownian motions. However, to prove  independence between them, we need to follow a different approach. We start by writing
$$\Delta(u,v)-\Delta(0,0)=\Gamma(u)+\Gamma^u(v)\,,$$
where
$$\Gamma(u):=\left(\Delta(u,0)-\Delta(0,0)\right)\,\mbox{ and }\,\Gamma^u(v):=\left(\Delta(u,v)-\Delta(u,0)\right)\,.$$
We will prove that $\epsilon^{-1/2}\Gamma(\epsilon u)$ and $\epsilon^{-1/2}\Gamma^{\epsilon u}(\epsilon v)$ converges jointly to a pair of independent standard Brownian motions, in the topology of continuous functions on compact sets. Since we already have tightness, we only need to prove finite dimensional convergence.
\newline

For $\rho\in(0,1)$, let $\s^1_\rho$ and $\s^2_\rho$ be two independent copies of the stationary profile \eqref{defstat}. Set $m=\lfloor 2^{-1}n\rfloor$ and define
$$C_m^{1,[k]_n}:=\left\{j\in\ZZ\,:\, (m+j,m-j)\leq [k]_n\right\}\,,$$
and 
$$C_m^{2,[k]_0}:=\left\{j\in\ZZ\,:\, (m+j,m-j)\geq [k]_0\right\}\,.$$
Define also 
$$L_m^{\rho,1}[k]_n:=\max_{j\in C^{[k]_n}_1} \left\{ \s^1_\rho(j) + L\left([j]_m,[k]_n\right)\right\}\,,$$
and 
$$L_m^{\rho,2}[k]_0:=\max_{j\in C^{[k]_0}_2} \left\{ \s^2_\rho(j) + L\left([k]_0,[j]_m\right)\right\}\,.$$
In words, $L_m^{\rho,1}$ is a last-passage percolation model with initial profile $\s^1_\rho$ (along $x+y=2m$), while $L_m^{\rho,2}$ is a reversed last-passage percolation model with initial profile $\s^2_\rho$. Thus, $L_m^{\rho,1}$ and $L_m^{\rho,2}$ are independent. We will approximate the local fluctuations of $L_{un^{2/3}}[vn^{2/3}]_n-L_{un^{2/3}}[0]_n$ by the local fluctuations of $L_m^{\rho,1}[vn^{2/3}]_n-L_m^{\rho,1}[0]_n$, and the local fluctuations of $L_{un^{2/3}}[0]_n-L_{0}[0]_n$ by the local fluctuations of $L_m^{\rho,2}[un^{2/3}]_0-L_m^{\rho,2}[0]_0$, and then extend these approximations to $\Gamma^u(v)$ and $\Gamma(u)$.  
\newline

Let $Z^{j}_m[k]_n$ be such that   
$$[Z^{j}_m[k]_n]_m=(m+Z^{j}_m[k]_n,m-Z^{j}_m[k]_n)\,$$
is the intersection point between $\left\{(m+i,m-i)\,:\,i\in\ZZ\right\}$ and the geodesic from $[j]_0$ to $[k]_n$. We denote $Z_m^{\rho,1}[k]_n$ the exit point associated to $L^{\rho,1}[k]_n$.  
\begin{lem}\label{comparison1}
Let $k\leq l$ and $n\geq 1$. If $Z^{j}_m[l]_n\leq Z_m^{\rho,1}[k]_n$ then
$$L_{j}[l]_n-L_{j}[k]_n \leq L_m^{\rho,1}[l]_n-L_m^{\rho,1}[k]_n\,.$$
If $Z_m^{\rho,1}[l]_n\leq Z_m^{j}[k]_n$ then
$$L_m^{\rho,1}[l]_n-L_m^{\rho,1}[k]_n\leq L_{j}[l]_n-L_{j}[k]_n\,.$$ 
\end{lem}

\noindent{\bf Proof\,\,} We can write
$$L_j[l]_n-L_j[k]_n=(L_j[l]_n-L_j[0]_m)-(L_j[k]_n-L_j[0]_m)\,,$$
and since 
$$L_j[l]_n=L_j[Z^{j}_m[l]_n]_m+ L\left(Z^{j}_m[l]_n,[l]_n\right)\,,$$
we have that  
\begin{eqnarray*}
L_j[l]_n-L_j[0]_m &=& \left(L_j[Z^{j}_m[l]_n]_m-L_j[0]_m\right)+ L\left([Z^{j}_m[l]_n]_m,[l]_n\right)\\
&=& \rb_j\left(Z^{j}_m[l]_n\right)+L\left([Z^{j}_m[l]_n]_m,[l]_n\right)\,,
\end{eqnarray*}
where $\rb_j$ is the profile induced by $L_j$ along $x+y=2m$:
$$ \rb_j\left(z\right):= L_j[z]_m-L_j[0]_m \,.$$
Thus, $L_j[l]_n-L_j[0]_m$ can be formulated as in Lemma \ref{comparison} and the proof follows the same lines.

\hfill$\Box$\\ 

\begin{lem}\label{compcontrol2}
For $r\geq 0$ set 
$$\rho_n^{\pm}=\rho_n^{\pm}(r)=1/2\pm\frac{r}{n^{1/3}}\,,$$
and define the event
$$\bar E^1_{n,m}(r):=\left\{Z_m^{\rho_n^-,1}[n^{2/3}]_n\leq Z^{0}_m [0]_n\,\,\mbox{ and }Z_m^{n^{2/3}}[n^{2/3}]_n\leq Z_m^{\rho_n^+,1} [0]_n\,\, \right\}\,.$$
Then,  
$$\limsup_{r\to\infty}\limsup_{n\to\infty}\P\left(\bar E^1_{n,m}(r)^c\right)=0\,.$$
\end{lem}

\noindent{\bf Proof\,\,} By Theorem 2.5 \cite{BaCaSe}, 
$$\limsup_{r\to\infty}\limsup_{n\to\infty}\P\left(|Z_m^{j}[k]_n|\geq rn^{2/3} \right)=0\,,\,\mbox{for all $j,k\in [0,n^{2/3}]$ }\,,$$
which implies Assumption \ref{exitcontrol}. Therefore, Lemma \ref{compcontrol2} follows from Lemma \ref{compcontrol}. 

\hfill$\Box$\\ 

Define
$$\Gamma^u_n(v):=\frac{L_{un^{2/3}}[vn^{2/3}]_n-L_{un^{2/3}}[0]_n}{2^{3/2}}\,,$$
$$B^{\pm}_{n,1}(v):=\frac{L_m^{\rho_n^\pm,1}[vn^{2/3}]_n-L_m^{\rho_n^\pm,1}[0]_n-\mu_{\rho_n^\pm}vn^{2/3}}{2^{3/2}n^{1/3}}\,,$$ 
and $B^{\pm,\epsilon}_{n,1}(v):=\epsilon^{-1/2}B^{\pm,1}_{n,1}(\epsilon v)$. By Lemma \ref{comparison1}, on the event $\bar E^1_{n,m}(r)$ (analog to \eqref{localBM})
\begin{equation*}
B^{-,\epsilon}_{n,1}(v)-3\sqrt{2}A\epsilon^{1/2} r \leq \epsilon^{-1/2}\Gamma^{\epsilon u}_n(\epsilon v)\leq B^{+,\epsilon}_{n,1}(v)+3\sqrt{2}A\epsilon^{1/2} r\,,
\end{equation*}
for all $(u,v)\in[0,A]^2$ (and $\epsilon A\in [0,1]$). If $r_\epsilon:=\epsilon^{-\beta}$, with $\beta\in(0,1/2)$, and $(u_i,v_i)\in[0,A]^2$ then 
$$\P\left(\cap_{i=1}^j\left\{ \epsilon^{-1/2}\Gamma^{\epsilon u_i}(\epsilon v_i)\leq x_i\right\}\right)\leq \P\left(\cap_{i=1}^j\left\{ B_1(v_i)\leq x_i +3\sqrt{2}A\epsilon^{1/2-\beta}\right\}\right)+\limsup_{n\to\infty}\P\left(\bar E^1_{n,m}(r_\epsilon)^c\right)\,  $$
and
$$\P\left(\cap_{i=1}^j\left\{ \epsilon^{-1/2}\Gamma^{\epsilon u_i}(\epsilon v_i)\leq x_i\right\}\right)\geq \P\left(\cap_{i=1}^j\left\{ B_1(v_i)\leq x_i -3\sqrt{2}A\epsilon^{1/2-\beta}\right\}\right)+\limsup_{n\to\infty}\P\left(\bar E^1_{n,m}(r_\epsilon)^c\right)\,.$$
By Lemma \ref{compcontrol2}, this proves finite-dimensional convergence:
$$\lim_{\epsilon\downarrow 0}\P\left(\cap_{i=1}^j\left\{  \epsilon^{-1/2}\Gamma^{\epsilon u_i}(\epsilon v_i)\leq x_i\right\}\right)= \P\left(\cap_{i=1}^j\left\{ B_1(v_i)\leq x_i \right\}\right)\,.$$

To prove that 
$$\lim_{\epsilon\downarrow 0}\P\left(\cap_{i=1}^j\left\{  \epsilon^{-1/2}\Gamma(\epsilon u_i)\leq x_i\right\}\right)= \P\left(\cap_{i=1}^j\left\{ B_2(u_i)\leq x_i \right\}\right)\,,$$
where $B_2$ is a Brownian motion that is independent of $B_1$, one has to prove the analogs of Lemma \ref{comparison1} and Lemma \ref{compcontrol2}, and use the approximation   
$$B^{\pm}_{n,2}(u):=\frac{L_m^{\rho_n^\pm,2}[un^{2/3}]_0-L_m^{\rho_n^\pm,2}[0]_0-\mu_{\rho_n^\pm}vn^{2/3}}{2^{3/2}n^{1/3}}\,,$$
that is clearly independent of $B^{\pm}_{n,1}$.

\noindent\paragraph{\bf Acknowledgments.} The author would like to thank Eric Cator, Ivan Corwin and Jeremy Quastel for enlightening discussions about  last-passage percolation and Airy sheets, and Patrick Ferrari and Timo Sepp\"al\"ainen for useful comments and corrections of a previous version of this article. The author was partially supported by the CNPQ grants 421383/2016­0 and 302830/2016-2, and by the FAPERJ grant E-26/203.048/2016.

\end{document}